\documentclass[11pt]{article}
\usepackage{epsfig,amsmath,amsthm,amssymb,graphics,amsfonts,psfrag}
\usepackage[letterpaper,margin=1in]{geometry}

\theoremstyle{plain}

\theoremstyle{definition}

\def\R{\mathbb{R}}

\def\N{\mathbb{N}}
\def\I{\infty}

\newcommand{\be}{\begin{equation}}
\newcommand{\ee}{\end{equation}}
\newcommand{\bea}{\begin{eqnarray}}
\newcommand{\eea}{\end{eqnarray}}
\newcommand{\beann}{\begin{eqnarray*}}
\newcommand{\eeann}{\end{eqnarray*}}
\newcommand{\benn}{\begin{equation*}}
\newcommand{\eenn}{\end{equation*}}

\def\ra{\rightarrow}
\def\I{\infty}

\begin{document}
 
\title{On decomposing mixed-mode oscillations and their return maps}
\author{Christian Kuehn\thanks{Max Planck Institute for the Physics of Complex Systems}}

\maketitle

\begin{abstract}
Alternating patterns of small and large amplitude oscillations occur in a wide variety of physical, chemical, biological and engineering systems. These mixed-mode oscillations (MMOs) are often found in systems with multiple time scales. Previous differential equation modeling and analysis of MMOs has mainly focused on local mechanisms to explain the small oscillations. Numerical continuation studies reported different MMO patterns based on parameter variation. This paper aims at improving the link between local analysis and numerical simulation. Our starting point is a numerical study of a singular return map for the Koper model which is a prototypical example for MMOs that also relates to local normal form theory. We demonstrate that many MMO patterns can be understood geometrically by approximating the singular maps with affine and quadratic maps. Motivated by our numerical analysis we use abstract affine and quadratic return map models in combination with two local normal forms that generate small oscillations. Using this decomposition approach we can reproduce many classical MMO patterns and effectively decouple bifurcation parameters for local and global parts of the flow. The overall strategy we employ provides an alternative technique for understanding MMOs.
\end{abstract}


\textbf{Complex oscillatory patterns have been observed in a wide variety of applications. Analyzing these patterns from a dynamical perspective has been an active area of research for decades. However, several mathematical breakthroughs in the last 15 years have provided substantial additional insight into phenomena that describe local oscillations. In the present paper, we provide a numerical study of the singular Poincar\'{e} map in the Koper model. We demonstrate that many MMO patterns for the Koper model can already be understood just using approximations of singular limit maps. The results for the Koper model suggest that a local-global numerical simulation approach combining normal forms with discrete maps can be effective. We show that this abstract approach reproduces many typical MMO patterns that have been observed in applications. This methodology aims to close a gap between previous numerical studies of MMO patterns and analytical results about local normal forms.}

\section{Introduction}  

Mixed-mode oscillations (MMOs) are patterns of small and large amplitude oscillations in a time series that differ at least by one order in magnitude. They have been observed experimentally in the Belousov-Zhabotinsky reaction in the 1970's and 1980's \cite{HudsonHartMarinko,MaselkoSwinney} and have been encountered more recently in a wide variety of different experiments \cite{HauckSchneider,HauserOlsen,OrbanEpstein,Dicksonetal}. The basic classification has been based on counting the number of small oscillations $s$ and large oscillations $L$ so that we can symbolically represent an MMO by
\benn
\ldots L_{j-1}^{s_{j-1}}L_{j}^{s_{j}}L_{j+1}^{s_{j+1}}\ldots
\eenn
where $j\in \N$ is an index. For example, if we have a periodic time series that has 2 large amplitude oscillations (LAOs) and then 3 small amplitude oscillations (SAOs) we get $\ldots 2^32^32^3\ldots$ or simply an MMO of type $2^3$. Systems exhibiting MMOs are often modeled using differential equations \cite{DegnOlsenPerram,Barkley}. Local bifurcation theory \cite{BronsKrupaWechselberger,Guckenheimer7} and numerical methods \cite{Koper,DesrochesKrauskopfOsinga2} have been developed to gain a lot of insights into SAO generating mechanisms \cite{BronsKaperRotstein}. A multiple time scale structure of the system is a key component for many local mechanisms. A detailed survey of this theory and its applications to particular models has been completed recently \cite{Desrochesetal}. The main findings of many numerical studies (see e.g. \cite{Koper,WechselbergerWeckesser,RotsteinWechselbergerKopell,DesrochesKrauskopfOsinga1}) and experiments (see e.g. \cite{HudsonHartMarinko,MaselkoSwinney1,VenrooijKoper}) are transition sequences of periodic orbits; for example, if we only consider MMOs with patterns of the form $\cdots L^sL^s\cdots$ such a transition sequence can be represented as follows
\be
\label{eq:MMO_patt_main}
\cdots \ra (L_{p_1})^{s_{p_1}} \ra (L_{p_2})^{s_{p_2}} \ra (L_{p_3})^{s_{p_3}} \ra \cdots
\ee
where $p$ is a control/bifurcation parameter i.e. under variation of a single parameter changing patterns of MMOs can be observed. To understand patterns of the form \eqref{eq:MMO_patt_main} several approaches have been used. The theory of local normal forms has been applied to explain the SAOs and then it is usually assumed that the global return mechanism satisfies certain properties (see e.g. \cite{WechselbergerWeckesser,BronsKrupaWechselberger}) so that the local theory becomes applicable or a phenomenological model for the return map is proposed \cite{MilikSzmolyan,MilikSzmolyanLoeffelmannGroeller}. Another approach is to compute Poincar{\'{e}} maps \cite{Kuznetsov} numerically under parameter variation (see e.g. \cite{KawczynskiKhavrusStrizhak,MedvedevYoo,GuckenheimerScheper}) to explain transitions of MMO patterns or to use numerical continuation \cite{KrauskopfOsingaGalan-Vioque} to subdivide parameter space (see e.g. \cite{Koper,DesrochesKrauskopfOsinga2}). These techniques have provided tremendous insight into what types of sequences \eqref{eq:MMO_patt_main} can be found in different systems. However, all previously mentioned studies vary parameters in such a way that local and global dynamics change \textit{simultaneously}. Here we suggest that to understand which patterns of the form \eqref{eq:MMO_patt_main} occur one also has to ask what happens when this parameter coupling is not present. Only in this context one is able to distinguish the effects of parameter variation on the local normal form from the variation of parameters in the Poincar\'{e} map. We start by applying this idea in the context of Koper's model \cite{Koper}. For Koper's model the local dynamics is well-understood \cite{Desrochesetal} and SAOs are generated by folded nodes \cite{SzmolyanWechselberger1,Wechselberger2} and folded saddle-nodes of type II (or singular Hopf bifurcation, \cite{Guckenheimer7,Desrochesetal}) which are normal forms for systems fast-slow systems with three variables (see also Appendices \ref{ssec:folded_nodes}-\ref{ssec:sing_Hopf} for a brief review).\\

\textit{Remark:} We point out that folded nodes and singular Hopf bifurcation are two possible normal forms under the assumptions of fast-slow systems structure and non-degeneracy assumptions for a folded critical manifold. Obviously one can also suggest other possible SAO mechanisms \cite{GuckenheimerWillms,MedvedevYoo}. However, we have chosen to focus on the Koper model that is well-described locally by the two normal forms described above. The main reasons for this choice are that many experimental and analytical studies have been found that exhibit folded nodes and/or singular Hopf bifurcation (see the review \cite{Desrochesetal} for a list systems with folded nodes and singular Hopf bifurcation). Furthermore, it has recently been shown that both mechanisms also relate to delayed Hopf bifurcation \cite{KrupaWechselberger} which has been proposed as another SAO mechanism.\\ 

The global return mechanism for MMOs in the Koper model is provided by a cubic relaxation-oscillation mechanism \cite{MKKR,SzmolyanWechselberger} that has already been investigated by van der Pol in the 1920s \cite{vanderPol,vanderPol1}. Here we provide numerical computations of the global Poincar\'{e} return map as a composition of several maps in the singular limit of perfect time scale separation. These calculations reveal that the return map can be surprisingly regular. Using affine and quadratic approximations to the singular maps we investigate MMO patterns and find that the approximations suffice to understand MMO sequences observed in extensive numerical continuation. Motivated by these results we combine two local normal form ODEs with abstract linear and quadratic maps to study MMOs. It is shown that classical sequences of the form \eqref{eq:MMO_patt_main} as well as chaotic MMOs can be easily generated in this framework. In particular, it is easy to design MMO patterns and to understand the differences in local and global parameter effects. We point out that this study also contributes to closing the gap between numerical simulation and local normal forms by reproducing several of the MMO transition sequences observed by a simultaneous local and global parameter variation in the Koper model.\\

The paper is structured as follows. Appendix \ref{ssec:review} contains the necessary background for readers not familiar with fast-slow system and MMO generating mechanisms in these systems. The main part of this paper starts in Section \ref{sec:Koper} where the Koper model is introduced and its basic properties are reviewed. In Section \ref{sec:maps1} the global singular return map for the Koper model is decomposed into several more tractable flow maps using numerical simulations. In Section \ref{sec:maps2} the maps are approximated by affine and quadratic map models; Appendix \ref{sec:error} contains a discussion of the approximation error. In Section \ref{sec:MMO} the global aspects of MMOs in the Koper model are analyzed using the flow map models. In Section \ref{sec:loc_glob} we consider a standard local-global decomposition of the MMO generating mechanisms. The key point is that we suggest to separate the parameter dependencies for the local and global models. We combine a global return map model with local SAOs induced by folded node and singular Hopf normal forms. We conclude with a brief outlook, describing the wider applicability of our approach, in Section \ref{sec:outlook}.

\section{The Koper Model}
\label{sec:Koper}

One version of the Koper model for MMOs is given by
\bea
\label{eq:Koper}
\epsilon_1 \dot{x} &=& y-x^3+3x, \nonumber\\
\dot{y} &=& kx-2(y+\lambda)+z,\\
\dot{z}&=& \epsilon_2(\lambda+y-z),\nonumber
\eea
where $(k,\lambda)$ are the main bifurcation parameters and $(\epsilon_1,\epsilon_2)$ are the singular perturbation parameters. The equations were first studied as a two-dimensional model by Boissonade and De Kepper \cite{BoissonadeDeKepper} modeling a prototypical chemical reaction. Koper \cite{Koper} added a third variable to a planar system and used numerical continuation techniques \cite{Govaerts,Doedel_AUTO2007} to study MMOs \cite{Koper}. It is very important to note that equations similar or equivalent to \eqref{eq:Koper} have been proposed many times \textit{independently} by several different research groups \cite{GoryachevStrizhakKapral,StrizhakKawczynski,KawczynskiStrizhak,BronsKrupaWechselberger,KrupaPopovicKopell,Guckenheimer7}
as a ``canonical'', ``minimal'' or ``typical'' model for MMOs. The version \eqref{eq:Koper} of Koper's model was proposed by the author and co-workers in \cite{Desrochesetal}; it is obtained by a coordinate transformation of Koper's original model and has the symmetry
\benn
(x,y,z,\lambda,k)\mapsto (-x,-y-z,-\lambda,k)
\eenn
which allows us to restrict to parameter regions with $\lambda\geq0$ or $\lambda\leq 0$ without loss of generality. Here we shall only review the local bifurcation structure briefly an introduce the necessary notation; a detailed local fast-slow systems analysis of \eqref{eq:Koper} can be found in \cite{Desrochesetal}. We also point out that the terminology reviewed in Appendix \ref{ssec:review} will be assumed from now on.

\begin{figure}[htbp]
	\centering
		\includegraphics[width=0.8\textwidth]{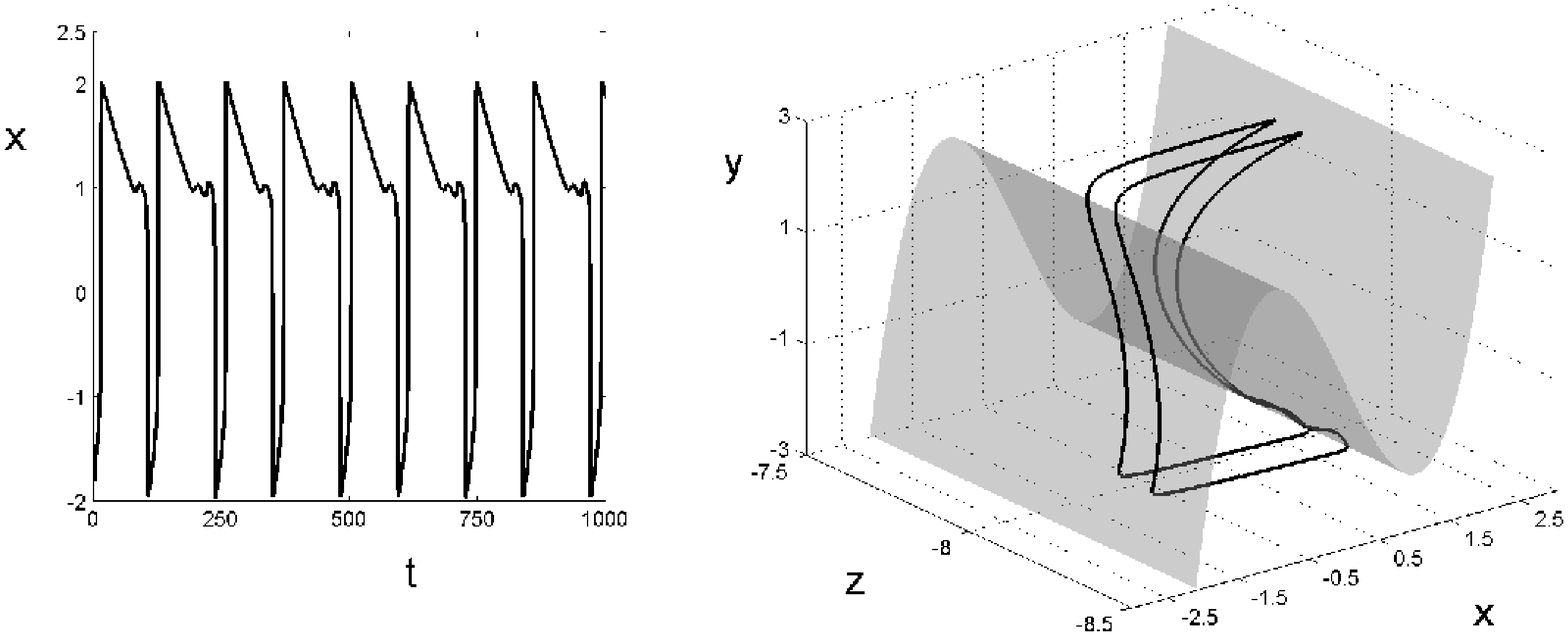}
	\caption{\label{fig:fig9}The parameter values for the simulation are $(\epsilon,k,\lambda)=(0.01,-10,-7)$. The time series for a $1^11^2$ MMO for the variable $x$ is shown on the left and the phase space trajectory is shown on the right; the critical manifold $C_0$ is shown in grey.}	
\end{figure}

If $0<\epsilon_{1,2}\ll 1$ holds then \eqref{eq:Koper} is a three time-scale system. We shall focus on the case $\epsilon_2=1$ and $0<\epsilon_1=:\epsilon\ll 1$ in which case we have one fast variable $x$ and two slow variables $(y,z)$. The critical manifold is
\benn
C_0=\{(x,y,z)\in\R^3:y=x^3-3x=:c(x)\}.
\eenn
The typical cubic (or S-shaped) structure splits the critical manifold into several parts
\benn
C_0=C^{a,-}\cup F_-\cup C^r\cup F_+\cup C^{a,+}
\eenn 
where $C^{a,-}:=C\cap \{x<-1\}$, $C^{a,+}:=C\cap \{x>1\}$ are normally hyperbolic attracting, $C^r:=C\cap\{-1<x<1\}$ is normally hyperbolic repelling and 
\benn
F_-:=C\cap \{x=-1\}=\{(-1,2,z)\}\qquad \text{and} \qquad F_+:=C\cap \{x=1\}=\{(1,-2,z)\}
\eenn 
are fold curves curves of the critical manifold. MMOs can easily be observed in simulations; see Figure \ref{fig:fig9}. The desingularized slow subsystem is
\bea
\label{eq:desing_sf}
\dot{x}&=& kx-2(c(x)+\lambda)+z,\nonumber\\
\dot{z}&=& (3x^2-3)(\lambda+c(x)-z).
\eea

There are two folded singularities
\benn
p_\pm=(\pm 1,2\lambda \mp(4+k)).
\eenn 
By symmetry we shall only focus on the folded singularity $p_+$. The linearization of the desingularized slow flow $p_+$ is
\be
\label{eq:Jacobian1}
\left(\begin{array}{c}\dot{X} \\ \dot{Z} \\ \end{array}\right)= \left(\begin{array}{cc}
-10 & 1 \\
-6(8 + \lambda) & 0\\
\end{array} \right)\left(\begin{array}{c}X \\ Z \\ \end{array}\right) =: A_+ \left(\begin{array}{c}X \\ Z \\ \end{array}\right)
\ee 
where we already set $k=-10$ which will fixed from now on. Note that in this case $\{X=0 \}$ corresponds to the fold line $F_+$ and $p_+$ is located at the origin. The eigenvalues of $A_+$ are
\benn
\sigma_w(\lambda)=-5+\sqrt{-23-6\lambda} \qquad \text{and} \qquad \sigma_s(\lambda)=-5-\sqrt{-23-6\lambda},
\eenn
with associated eigenvectors 
\benn
\Sigma_w(\lambda)=\left(\begin{array}{c}\frac{1}{5+\sqrt{-23-6\lambda}} \\ 1 \\ \end{array}\right) \qquad \text{and} \qquad
\Sigma_s(\lambda)=\left(\begin{array}{c}\frac{1}{5-\sqrt{-23-6\lambda}} \\ 1 \\ \end{array}\right). 
\eenn
Therefore $p_+$ is a folded saddle for $\lambda<-8$, a folded saddle-node of type II (FSN II \cite{SzmolyanWechselberger1,Desrochesetal}) for $\lambda_{\text{FSN II}}=-8$ and a folded node for $\lambda\in(-8,-23/6)$. At $\lambda_{nf}=-23/6$ the transition from a folded node to a folded focus occurs. The singular Hopf bifurcation for the full system occurs $O(\epsilon)$ away from $\lambda_{\text{FSN II}}$. It is supercritical and the stable global equilibrium $q$ loses stability at this point. Therefore the interesting parameter region for MMOs is 
\benn
\lambda\in(\lambda_{\text{FSN II}},\lambda_{nf})=(-8,-23/6),\qquad k=-10.
\eenn
Note that this parameter region represents a typical one-parameter MMO sequence \cite{Koper,Desrochesetal}. The important eigenvector for global returns is $\Sigma_s$ associated to the strong primary canard $\gamma_s$ as it bounds the rotational sectors lying on $C^{a,+}$. The $x$-component $\Sigma_s^x$ of $\Sigma_s$ lies, for the scaling we have chosen, between $\Sigma_s^x(-8)=\I$ and $\Sigma_s^x(-23/6)=\frac15$. Therefore the rotational sectors \cite{BronsKrupaWechselberger} that subdivide the funnel region are given by a convex cone with opening angle between $\frac{\pi}{2}$ and $\cos^{-1}(25/26)$.

\section{Return Maps - Decomposition}
\label{sec:maps1}

\begin{figure}[htbp]
	\centering
		\includegraphics[width=0.9\textwidth]{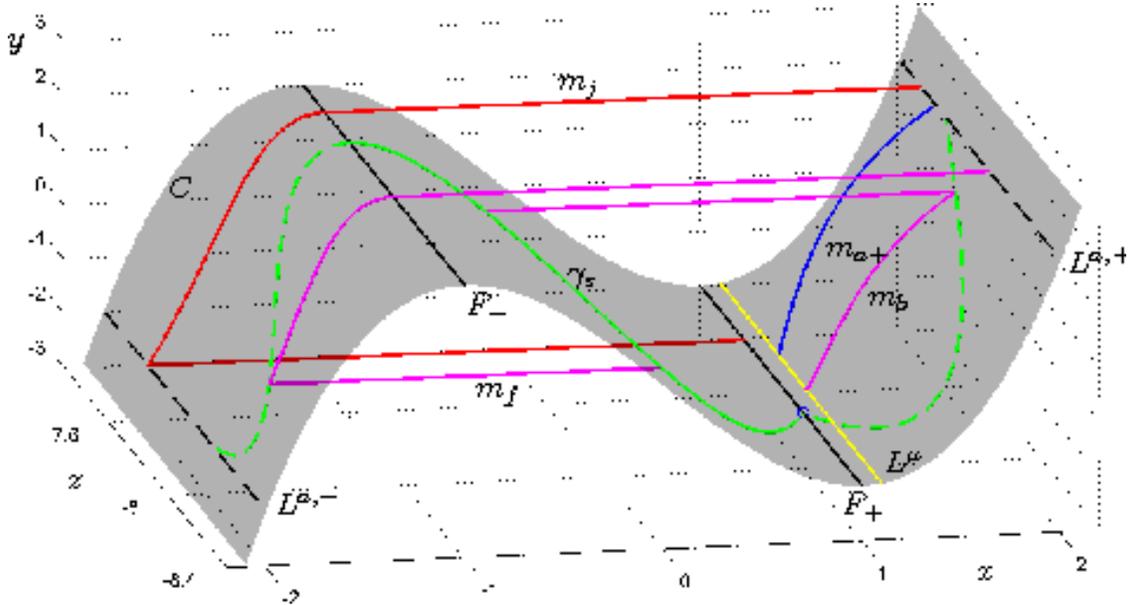}
	\caption{\label{fig:fig1}Illustration of the singular map decomposition; parameter values are $(\epsilon,k,\lambda)=(0,-10,-7)$. Definitions of all maps and domains are given at the beginning of Section \ref{sec:maps1}. Here we show: the critical manifold $C$ (grey), the strong canard $\gamma_s\subset C^m$ (green) and its projections to $C^{a,\pm}$ (dashed green), the fold lines $F_\pm$ (black) and their projections $L^{a,\mp}$ (dashed black, $\mu=0.1$), the line $L^\mu$ (yellow) and the folded node $p_+$ (blue circle). Examples for the maps $m_j$ (red, regular jump), $m_{a,+}$ (blue, flow towards $F_+$), $m_f$ (magenta, jump forward canard) and $m_b$ (magenta, jump backward canard) are displayed as well.}	
\end{figure}

Our goal is to analyze the structure of the global singular return map. Instead of using the standard approach of computing the Poincar\'{e} map between two fixed sections \cite{Shilnikov,GH,Kuznetsov} we are going to decompose the map according to fast-slow systems theory (see e.g. \cite{Haiduc1,GuckenheimerWechselbergerYoung,SzmolyanWechselberger} for this approach). We fix $\epsilon=0$ and recall that $k=-10$. Then we focus on $\lambda$ as the primary bifurcation parameter. In this case the folded node $p_+$ and the unique equilibrium $q$ account for the SAOs. Observe that global returns to a neighborhood of $p_+$ can be decomposed. See Figure \ref{fig:fig1} for an illustration of the one-dimensional singular maps we are going to define:

\begin{enumerate}
 \item[(a)] Trajectories can reach the fold line $F_+$ at a jump point and follow the fast flow to the drop curve $L^{a,-}:=C\cap \{x=-2\}$. Then trajectories follow the slow flow induced by \eqref{eq:desing_sf} to $F_-$ and jump to the drop curve $L^{a,+}:=C\cap\{x=2\}$. We denote this map by
\benn
m_j:F_+\ra L^{a,-}\ra F_-\ra L^{a,+}
\eenn
where $j$ indicates that we consider a regular jump. We denote the intermediate map by $m_{a-}:L^{a,-}\ra F_-$. Observe that if we parametrize the domain and range by $z$ then the intermediate map $m_{a-}$ is the only non-trivial component of the map $m_j$ and the other parts of $m_j$ are the identity with respect to $z$. 
 \item[(b)] Trajectories can flow into the folded node $p_+$. Suppose we consider trajectories tracking the part of the strong canard $\gamma_s$ contained in $C^r$. These trajectories jump at some point from $\gamma_s$ to $C^{a,-}$ and flow into $F_-$ before jumping to $L^{a,+}$. Denote this map by
\benn
m_f:\gamma_s\ra C^{a,-}\ra F_-\ra L^{a,+} 
\eenn 
where $f$ indicates a jump forward (or away) singular canard orbit; again observe that only the part $C^{a,-}\ra F_-$ is non-trivial with respect to $z$.
\item[(c)] Trajectories tracking the strong canard $\gamma_s\subset C^r$ can also jump at some point from $\gamma_s$ to $C^{a,+}$ and flow into $F_+$. It will be advantageous to terminate this map at a line $L^\mu:=C\cap \{x=1+\mu\}$ for some $\mu\geq0$ sufficiently small. Then we have a map 
\benn
m_b:\gamma_s\ra C^{a,+}\ra L^\mu
\eenn  
where $b$ indicates a jump backwards (or back) singular canard orbit.
\item[(d)] There is also a map induced by the slow flow on $C^{a,+}$ starting from the drop curve $L^{a,+}$ towards the fold line 
\benn
m_{a,+}:L^{a,+}\ra L^\mu
\eenn
\item[(e)] The linearization \eqref{eq:Jacobian1} at the folded singularity $p_+$ can be used to define a flow map in the fold region
\benn
m_s:L^\mu\ra F_+
\eenn
\end{enumerate}

\begin{figure}[htbp]
	\centering
		\includegraphics[width=0.95\textwidth]{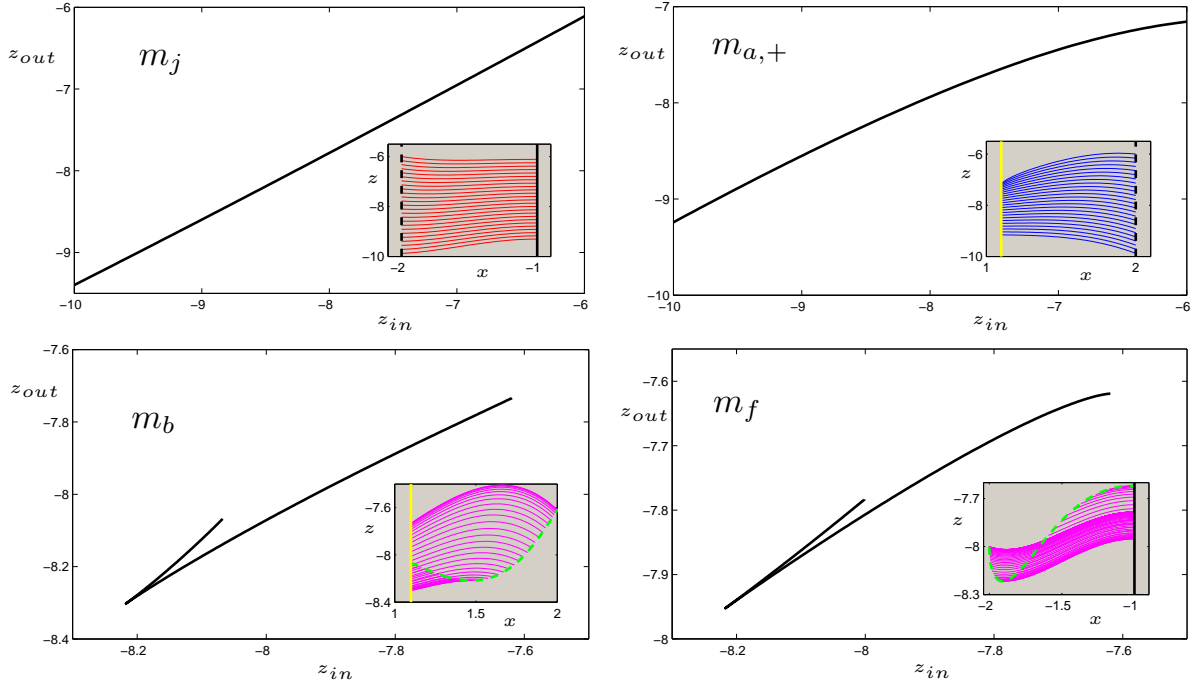}
	\caption{\label{fig:fig2}Singular maps for $(k,\lambda)=(-10,-7)$ with respect to the $z$-variable i.e. the horizontal axis shows $z=z_{in}$ and the vertical axis shows $z_{out}=m_{K}(z_{in})$ for $K\in\{j,(a,+),b,f\}$. The insets (grey background) illustrate the phase space flow on the attracting critical manifolds $C^{a,\pm}$ associated to the maps $m_K$; we show only every tenth trajectory in the computation of $m_K$. The color coding is the same as in Figure \ref{fig:fig1}.}	
\end{figure}

Figure \ref{fig:fig2} shows representatives of the maps $m_j$, $m_{a,+}$, $m_b$ and $m_f$ for $\lambda=-7$ with respect to the variable $z$ and also the associated slow flows. The main observation is that the maps are surprisingly regular. 

\section{Return Maps - Modeling}
\label{sec:maps2}

\begin{figure}[htbp]
	\centering
		\includegraphics[width=0.8\textwidth]{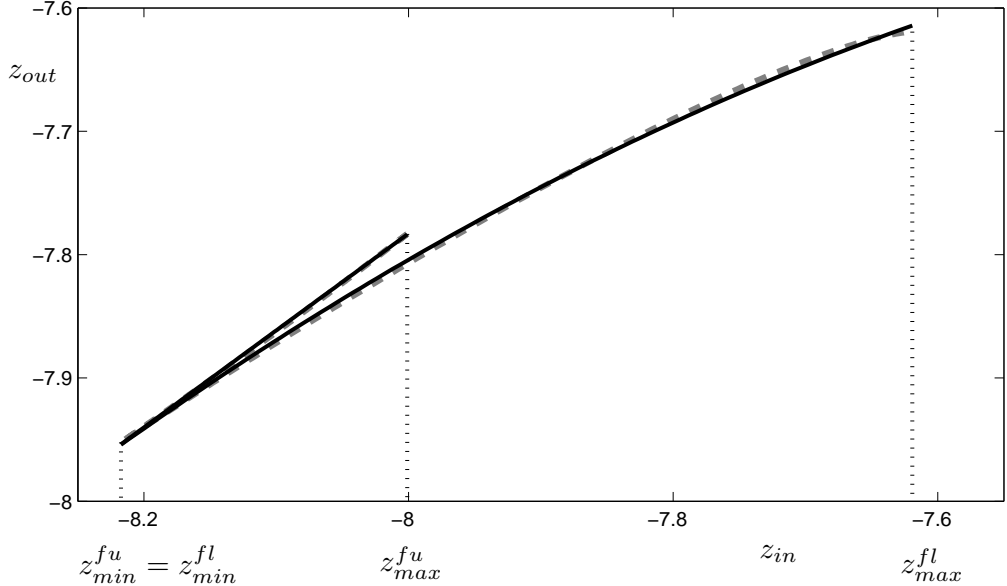}
	\caption{\label{fig:fig3}Singular map $m_f$ with $z_{out}=m_f(z_{in})$. The computed map is shown as a dashed grey curve and the approximations are shown in solid black (affine for upper part and quadratic for lower part). The bounds of the domains for each part of the map are marked as well (dotted vertical lines).}	
\end{figure}

In this section we are going to discuss the modeling of the maps computed in Figure \ref{fig:fig2}; a discussion of the approximation error as well as the error for $\epsilon>0$ is given in Appendix \ref{sec:error}. Figure \ref{fig:fig2} motivates considering affine and/or quadratic maps. The map $m_j$ seems to be close to an affine map which is due to the very simple regular slow flow from $L^{a,-}$ to $F_-$; see also Figure \ref{fig:fig1}. Similarly, we propose to model the map $m_{a,+}$ by a quadratic map. The maps induced from the projections of the strong canard $\gamma_s\subset C^r$ onto $C^{a,\pm}$ are multi-valued when parametrized with respect to $z$ due to the fold structure of $\gamma_s$; see Figure \ref{fig:fig1}. With another parametrization we expect that $m_b$ and $m_f$ are generically single-valued by uniqueness of solutions for the desingularized slow subsystem. The parametrization with respect to $z$ is very convenient. We propose to make the following ansatz:
\benn
m_K(z)=c_2(\lambda)z^2+c_1(\lambda)z+c_0(\lambda)
\eenn
for each map $m_K$ with $K\in\{j,(a,+),b,f\}$ where the coefficients $c_{0,1,2}(\lambda)$ are to be determined. We are going to illustrate the procedure for finding the coefficients for $m_f$ and just state the results we obtained for the other three maps. The ansatz is that $m_f$ can be decomposed as follows:
\be
\label{eq:map_mf}
m_f(z)=\left\{\begin{array}{ll} 
c^{fu}_1(\lambda) z + c^{fu}_0(\lambda) & \text{if $z^{fu}_{min}(\lambda) \leq z \leq z^{fu}_{max}(\lambda)$,}\\ 
c^{fl}_2(\lambda) z^2+ c^{fl}_1(\lambda) z + c^{fl}_0(\lambda) & \text{if $z^{fl}_{min}(\lambda) \leq z \leq z^{fl}_{max}(\lambda)$,}\\ \text{undefined} & \text{otherwise,}
 \end{array}\right.
\ee
where we impose continuity at the shared boundary point $m_f(z^{fu}_{min})=m_f(z^{fl}_{min})$. See Figure \ref{fig:fig3} for an example. In Figure \ref{fig:fig3} the upper part of $m_f$ is approximated by an affine map and the lower part by a quadratic. The computation of the approximation error in Appendix \ref{sec:error} for $\lambda\in(\lambda_{\text{FSN II}},\lambda_{nf})$ shows that for each fixed value of $\lambda$ the affine and quadratic models provide an approximation on the order of $10^{-2}$ of the singular maps obtained via numerical integration of slow flow trajectories on a fine mesh.\\  

\begin{figure}[htbp]
	\centering
		\includegraphics[width=0.90\textwidth]{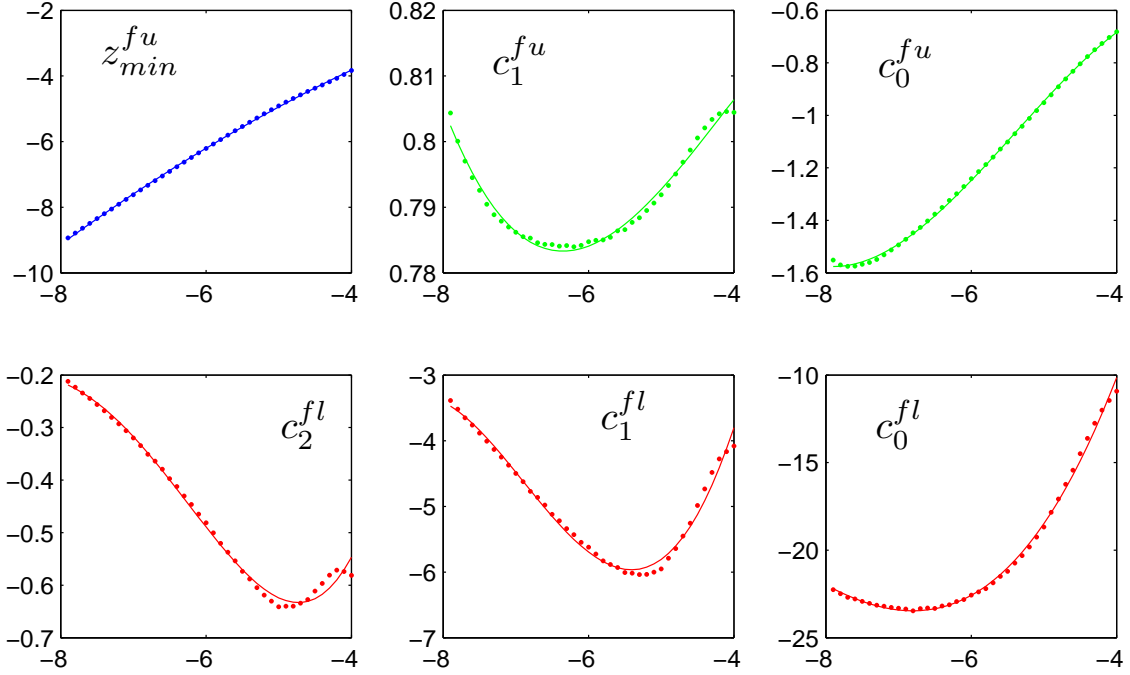}
	\caption{\label{fig:fig4}Horizontal axes are $\lambda$ and vertical axes are the respective coefficients e.g. top left figure shows $z^{fu}_{min}(\lambda)$. The only relevant part for the definition of \eqref{eq:map_mf} that is not shown is $z^{fl}_{max}$ which is as regular (almost linear) as the other parts of the domain boundaries for $m_f$. The dots are computed points and the curves provide polynomial fits (quadratic=blue, cubic=green and quartic=red).}	
\end{figure}

As a next step we investigate all functions depending on $\lambda$ in \eqref{eq:map_mf}. The boundary $z^{fu}_{max}$ is given by the folded singularity $p_+$ so that $z^{fu}_{max}(\lambda)=2\lambda+6$. We also know from the definition of \eqref{eq:map_mf} that $z^{fu}_{min}=z^{fl}_{min}$. The other functions of $\lambda$ can only be approximated numerically due to the nonlinear slow flows on $C^r$, which defines $\gamma_s$, and on $C^{a,-}$, which defines the map to $F_-$. Figure \ref{fig:fig4} shows numerical computations of the unknown functions of $\lambda$ in the definition of $m_f$ in \eqref{eq:map_mf}.  

Several observations can be made from Figure \ref{fig:fig4} and the previous remarks. All the domain boundaries $z^{fu}_{min}(\lambda)$, $z^{fu}_{max}(\lambda)$, $z^{fl}_{min}(\lambda)$ and $z^{fl}_{max}(\lambda)$ are very regular and seem to depend almost linearly on $\lambda$. The coefficients of the linear and quadratic polynomials have substantial nonlinear dependencies on $\lambda$ for the entire range $\lambda\in(\lambda_{\text{FSN II}},\lambda_{nf})$. This implies that although affine and quadratic maps can be very good approximations at fixed parameter values it will be more difficult to analyze the global return maps inducing MMOs as parameter-dependent families. For the other maps $m_b$, $m_j$ and $m_{a,+}$ we propose the following approximations:

\bea
m_b(z)&=&\left\{\begin{array}{ll} 
c^{bu}_2(\lambda) z^2+ c^{bu}_1(\lambda) z + c^{bu}_0(\lambda) & \text{if $z^{bu}_{min}(\lambda) \leq z \leq z^{bu}_{max}(\lambda)$,}\\ 
c^{bl}_1(\lambda) z + c^{bl}_0(\lambda) & \text{if $z^{bl}_{min}(\lambda) \leq z \leq z^{bl}_{max}(\lambda)$,}\\ \text{undefined} & \text{otherwise,}
 \end{array}\right. \label{eq:map_mb} \\
m_{j}(z)&=& c^{j}_1(\lambda) z + c^{j}_0(\lambda) \label{eq:map_mj},\\
m_{a,+}(z)&=& c^{a}_2(\lambda) z^2+ c^{a}_1(\lambda) z + c^{a}_0(\lambda), \label{eq:map_ma}
\eea
where we impose continuity at the shared boundary point for $m_b$ i.e. $m_b(z^{bu}_{min})=m_b(z^{bl}_{min})$. As a next step we are going to calculate the map for the linearized desingularized slow flow near $p_+$. The intersection of the eigendirection of $\Sigma_s$ with $L^\mu=\{x=1+\mu\}$ is easily calculated as
\benn
(1,2\lambda+6)^T+\left(\mu,\frac{\mu}{\Sigma_s^x(\lambda)}\right)^T=(1+\mu,2\lambda+6+\mu(5-\sqrt{-23-6\lambda}))=:(1+\mu,z^\mu(\lambda)).
\eenn  
Hence all trajectories that arrive at $L^\mu$ with $z\geq z^\mu(\lambda)$ will stay in the funnel and reach $p_+$ while trajectories for $z<z^\mu(\lambda)$ will first reach the fold line $F_-$ and jump to $L^{a,-}$. To see where on $F_+$ the last class of trajectories ends up we could just solve \eqref{eq:Jacobian1}. Note however that there exists an approximation for $z<z^\mu(\lambda)$ that just amounts to projecting $(\mu,Z(0))$ parallel to $\Sigma_s$ onto $\{X=0\}$ which is given by
\benn
(\mu,Z(0))\mapsto \left(0,Z(0)-\frac{\mu}{\Sigma_s^x}\right). 
\eenn
Therefore we get the local representation for the map $m_s$ in $(X,Z)$-coordinates
\benn
m^{loc}_s(Z)=\left\{\begin{array}{ll} 2\lambda+6 & \text{if $z \geq z^\mu(\lambda)$,}\\ 2\lambda+6 + Z-\frac{\mu}{\Sigma_s^x} & \text{if $z< z^\mu(\lambda)$.}\\   \end{array}\right.
\eenn
If $z$ is the coordinate obtained in original coordinates without linearization then
\benn
m_s(z)=\left\{\begin{array}{ll} 2\lambda+6 & \text{if $z \geq z^\mu(\lambda)$}\\ z-\frac{\mu}{\Sigma_s^x} & \text{if $z< z^\mu(\lambda)$}\\   \end{array}\right.
\eenn
where the error is $\mathcal{O}(\mu)$ as $\mu\ra 0$. With the different maps available we can proceed to analyze how they can be used to explain the global returns that generate LAOs.

\section{Mixed-Mode Oscillations}  
\label{sec:MMO}

\begin{figure}[htbp]
	\centering
		\includegraphics[width=0.98\textwidth]{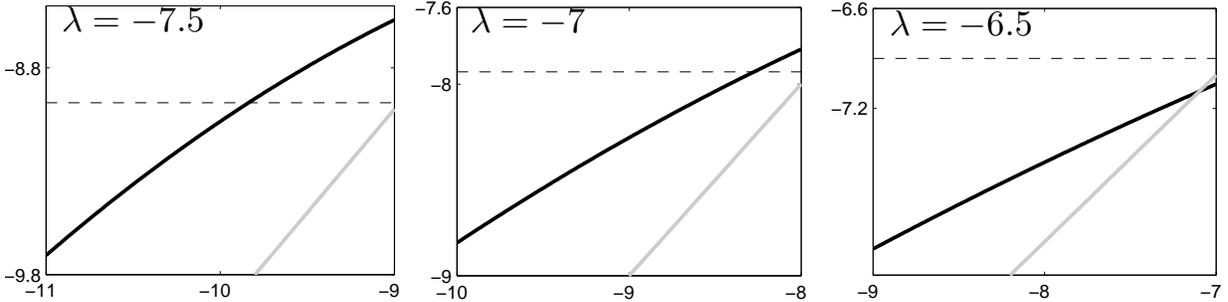}
	\caption{\label{fig:fig5}Map $(m_{a,+}\circ m_{j})(z)$, approximated by \eqref{eq:map_mj} and \eqref{eq:map_ma} with $k=-10$ and $\mu=0.1$. Horizontal axes are input $z$-coordinates on a domain $z\in((2\lambda+6)-2,2\lambda+6)\subset F_+$ and vertical axes are $(m_{a,+}\circ m_{j})(z)$ (think black curves). The location of the folded node funnel region $z^\mu(\lambda)$ is shown by horizontal dashed black lines and the diagonal is indicated by the thick grey line.}	
\end{figure}

Throughout this section we work with the polynomial approximations to the maps $m_{(.)}$ that have been derived in the last section. The first question we shall consider is what happens to trajectories that do not follow the canard $\gamma_s\subset C^r$ when arriving at $p_+$ or which land outside of the funnel region. The relevant map for this purpose is
\be
\label{eq:map_global1}
(m_{a,+}\circ m_{j}):F_+\cap \{z\leq 2\lambda+6\}\ra L^\mu
\ee
We are interested when part of the domain of \eqref{eq:map_global1} is returned inside the funnel so that $(m_{a,+}\circ m_{j})(z)>z^\mu(\lambda)$. Figure \ref{fig:fig5} shows the map \eqref{eq:map_global1} for three different values of $\lambda$. We observe that closer to the folded saddle-node of type II (i.e. near the singular Hopf bifurcation) trajectories that arrive outside the funnel on $F_+$ can get mapped back into the funnel under \eqref{eq:map_global1}. For $\lambda=-6.5$ in Figure \ref{fig:fig5} we observe that no trajectories can return into the funnel and that the return map $(m_s \circ m_{a,+}\circ m_{j})$ will have a stable fixed point since $\mu$ is small and hence the projection $m_s$ will preserve the intersection with the diagonal.\\

Hence we can consider several quantitative questions:

\begin{enumerate}
 \item For what values of $\lambda$ do trajectories from outside the funnel re-enter it?
 \item When does the map $(m_s \circ m_{a,+}\circ m_{j})$ have fixed points? When does the fixed point coincide 
with the folded node $p_+$?
 \item How are trajectories mapped into the funnel? More precisely, what is the dependence of the distance $\delta$ to the strong singular canard $\gamma_s\cap C^{a,+}$ upon varying $\lambda$?
\end{enumerate}

A trajectory starting for $z<2\lambda+6$ will re-enter the funnel after one global return if and only if
\beann
(m_{a,+} \circ m_j)(z)&=&c_2^a(\lambda)(c_1^j(\lambda)z+c_0^j(\lambda))^2+c_1^a(\lambda) (c_1^j(\lambda)z+c_0^j(\lambda))+c^a_0(\lambda)\\
&=& c_2^a (c_1^j)^2z^2+\left(2c_2^a c_1^j+c_1^a c_1^j\right)z+(c_0^j)^2c_2^a+c_1^ac_0^j+c_0^a<z^\mu(\lambda)
\eeann
By monotonicity of \eqref{eq:map_global1} on the required interval (see Figure \ref{fig:fig5}) we can just pick the folded node $z=2\lambda+6$ and determine when the condition fails; this yields the critical parameter value at which not all trajectories near $p_+$ return to the funnel in one iteration. We find that the parameter value at which $p_+$ gets returned to the boundary of the funnel is $\lambda=\lambda_r\approx -6.7887$. Next, we consider the fixed points of $(m_s \circ m_{a,+}\circ m_{j})$. Those points correspond to candidates representing relaxation oscillations. We find that at $\lambda=\lambda_r$ a stable fixed point appears for the map $(m_s \circ m_{a,+}\circ m_{j})$. Therefore we find that a transition to relaxation oscillations occurs near $\lambda_r$ for the full system and $\epsilon$ sufficiently small; this can be confirmed by numerical continuation \cite{Desrochesetal}. Note that the bifurcation that creates the fixed point occurs at the boundary of the domain of $(m_s \circ m_{a,+}\circ m_{j})$.\\

\begin{figure}[htbp]
	\centering
		\includegraphics[width=1\textwidth]{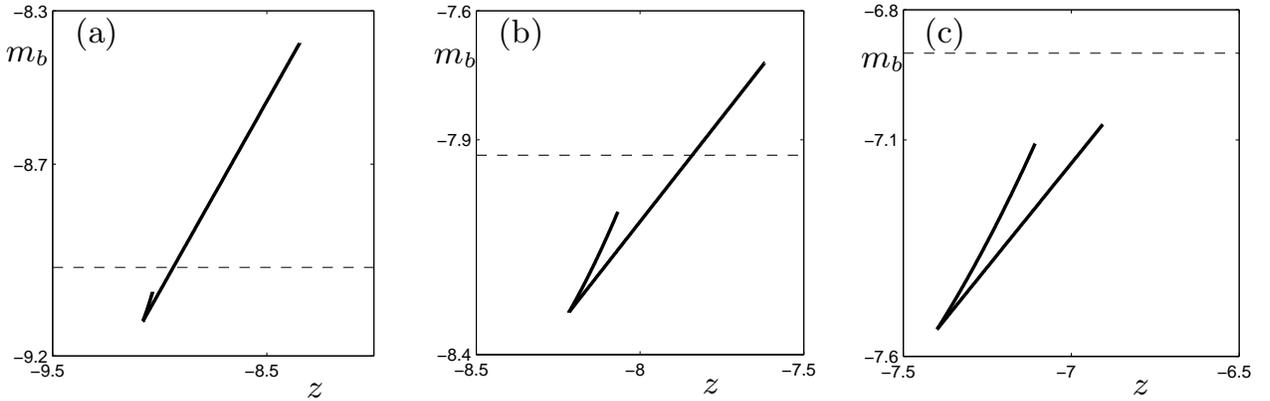}
	\caption{\label{fig:fig6}Singular maps with respect to the $z$-variable i.e. the horizontal axis shows $z=z_{in}$ and the vertical axis shows $m_b=m_{b}(z_{in})$. The horizontal dashed line indicates the funnel boundary $z^\mu(\lambda)$; here $\mu=0.1$. (a) $\lambda=-7.5$, (b) $\lambda=-7$ and (c) $\lambda=-6.5$.} 	
\end{figure}

As a next step we consider candidates that follow the canard $\gamma_s\cap C^r$ i.e. we consider the maps $m_f$ and $m_b$. We start with $m_b$ which represents medium-size canard-induced oscillations if trajectories from the domain of $m_b$ re-enter the funnel after one iteration step. Figure \ref{fig:fig6} plots three examples of the map $m_b$. The closer the parameter values are to the folded saddle-node of type II at $\lambda=-8$ the larger is the part of $\gamma_s\cap C^r$ that returns inside the funnel. The closer we are to relaxation oscillation at $\lambda=\lambda_r$ the more of $\gamma_s\cap C^r$ gets mapped outside the funnel.\\ 

Note that near $\lambda=-8$ with $\lambda>-8$ we must always have some part of $\gamma_s\cap C^r$ near $p_+$ that does get mapped outside the funnel since the opening cone angle of the funnel region is less than $\frac{\pi}{2}$; see Section \ref{sec:Koper}. Therefore there is always one part inside and one part outside the funnel for jump back canard orbits. Orbits in the full system that follow $\gamma^\epsilon_s$ for an $O(1)$-time on the slow time scale and get mapped back to $C^{a,+}$ via perturbation of $m_b$ represent intermediate oscillations.\\ 

\begin{figure}[htbp]
	\centering
		\includegraphics[width=1\textwidth]{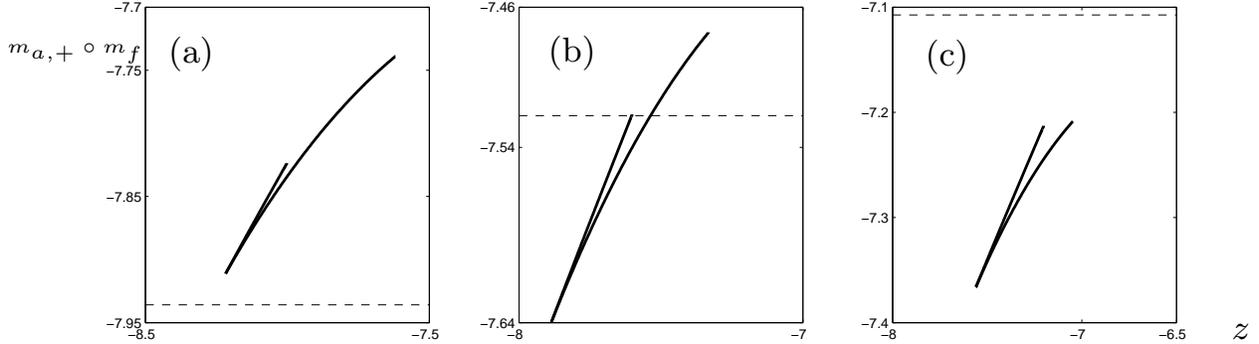}
	\caption{\label{fig:fig7}Singular maps with respect to the $z$-variable i.e. the horizontal axis shows $z=z_{in}$ and the vertical axis shows $m_{a,+}\circ m_f=(m_{a,+}\circ m_f)(z_{in})$. The horizontal dashed line indicates the funnel boundary $z^\mu(\lambda)$; here $\mu=0.1$. (a) $\lambda=-7$, (b) $\lambda=-6.8$ and (c) $\lambda=-6.6$.} 	
\end{figure}

For the map $m_f$ we immediately consider $m_{a,+}\circ m_f$ to see how jump forward canards get returned relative to the funnel. Figure \ref{fig:fig7} shows that there is a very rapid transition from jump forward canards that end all in the funnel for $\lambda=-7$ (Figure \ref{fig:fig7}(a)), a splitting of jump forward canards with respect to the funnel (Figure \ref{fig:fig7}(b)) and all jump forward canards outside the funnel for $\lambda=-6.6$ (see Figure \ref{fig:fig7}(c)). Let us consider the case when the entire jump forward canards end up in the funnel. This can be interpreted as a global MMO generating mechanism via canards. More precisely, a trajectory of the full system can make small oscillations near a folded node, follow $C^r\cap\gamma_s$ closely producing an intermediate oscillation and then return into the funnel. This provides a mechanism to transition small loops into large ones via canards. The closer we get to $\lambda=\lambda_r$ the more excursions outside the funnel occur which means that in this region we expect more mixed behavior of MMOs of type $L^s$ with $L>1$. It is also expected that period-doubling bifurcations of the return map can explain transitions between regions of different LAOs. Since resonances for the eigenvalues \cite{Wechselberger2} of the folded node are fewer near $\lambda_r$ we also expect $s$ to decrease if we increase $\lambda$. Hence we find that MMO sequences near a singular Hopf bifurcation will produce patterns with $s\gg 1$ and small $L$ while away from the singular Hopf $L^s$ patterns with $L\sim s$ are more likely to occur. All these findings agree with numerical continuation results in \cite{Koper,Desrochesetal}. 

Therefore one main conclusion from the numerical simulations considered here is that the singular limit decomposition is already sufficient to explain many MMO transition sequences. Indeed, in the singular limit we could already identify the local normal forms (see Section \ref{sec:Koper} and \cite{Desrochesetal}) and here we calculated a decomposition of the global return map. The main point is that we have used a different, and easily implementable, numerical \textit{technique} to understand geometrically many of the MMO patterns that have been found using extensive numerical continuation runs \cite{Koper}. 

\section{A Local-Global Model}
\label{sec:loc_glob}

We have seen that the global singular return maps for the Koper model are very regular and can often be described as affine or quadratic maps. The only feature of the global returns that is complicated to describe are canard orbits that follow the strong canard $\gamma_s\cap C^r$. These orbits describe intermediate oscillations i.e. orbits that, under parameter variation will grow to a large relaxation loop or decay to a small oscillation. However, many MMO transitions can be understood without these orbits as shown in the previous section. Hence it is natural to ask what happens if we do not consider these intermediate orbits and look at a simulation model for MMOs containing local and global maps. The local description of this model is chosen as a flow map for a folded node or a folded-saddle node ODE normal form; see Appendix \ref{ssec:review}. We assume without loss of generality that the folded singularity is located at the origin $(x,y,z)=(0,0,0)$. For the local dynamics we use the normal forms \eqref{eq:BKW} and \eqref{eq:Guck_SH}. Recall that the critical manifold of both normal forms is
\benn
C_0=\{(x,y,z)\in\R^3:y=x^2\}
\eenn
It is attracting for $x>0$ and repelling for $x<0$ and we denote the two branches of $C_0$ by $C^a_0$ and $C^r_0$. The associated attracting slow manifold provided by Fenichel Theory is
\benn
C^a_\epsilon=\{(x,y,z)\in \R^3:x=h^a_\epsilon(y,z)\}
\eenn
where the map $h^a_\epsilon$ is given by the implicit function theorem and $h^a_0(y,z)=\sqrt{y}$. Define two sections
\beann
\Sigma_1&:=&\{(x,y,z)\in\R^3|x=k_1\}\\
\Sigma_2&:=&\{(x,y,z)\in\R^3|x=-k_2\}
\eeann
for suitable fixed $k_j>0$, $k_j=O(\sqrt\epsilon)$ with $j=1,2$. The choice of scaling $O(\sqrt\epsilon)$ is prescribed by the fact that outside of a neighborhood of size $O(\sqrt\epsilon)$ of the origin Fenichel Theory applies. Define a map 
\be
\label{eq:loc_map_nf}
m_{12}:\Sigma_1\ra \Sigma_2
\ee
by the flow map of \eqref{eq:BKW} or \eqref{eq:Guck_SH}. Note that the sections $\Sigma_j$ are naturally parametrized by the coordinates $(y,z)$. The global return map $m_{21}:\Sigma_2\ra \Sigma_1$ will be modeled as follows:
\bea
m_{21}(y,z)&=&\left(\begin{array}{c} k_1^2 \\ m(z) \end{array}\right)+\epsilon\left[ \left(\begin{array}{cc} a_{11} & a_{12} \\ a_{21} & a_{22} \\  \end{array}\right) \left(\begin{array}{c}y \\ z \end{array}\right) + \left(\begin{array}{c} b_1 \\ b_2 \end{array}\right)\right]+O(\epsilon^2) \nonumber \\
&=&(k_1^2,m(z))^T+\epsilon [A(y,z)^T+b]+O(\epsilon^2) \label{eq:map_model_g}
\eea
where $m(z)=m_2z^2+m_1z +m_0$ and we require that the matrix $A$ is invertible. Note that we can make several further choices e.g. we could decide to include higher-order terms or to assume that $m(z)$ is modeled as an affine map and set $m_2=0$. The map $m_{21}$ has to satisfy a further constraint if we assume that all trajectories approach the origin exponentially close to the slow manifold $C^a_\epsilon$; this requires
\benn
k_2\sqrt\epsilon=(h^a_\epsilon\circ m_{21})(y,z).
\eenn
Another constraint to generate MMOs is that the global map $m_{21}$ maps some part of its domain close to the perturbation of the folded node funnel region. Although the model \eqref{eq:map_model_g} has formally nine free parameters $m_i$, $a_{jk}$, $b_l$ we can also view $m_{21}$ as an $O(\epsilon)$-perturbation of the leading order term which has only two or three parameters depending on the choice of model ($m_2=0$ or $m_2\neq 0$). Hence the description is low-dimensional, explicit and \textit{decouples} the local and global bifurcation structure of the problem. To illustrate the effect of global bifurcation parameters we numerically investigate two typical MMO sequences, one for each local normal form with fixed local parameters.\\

\begin{figure}[htbp]
	\centering
		\includegraphics[width=1\textwidth]{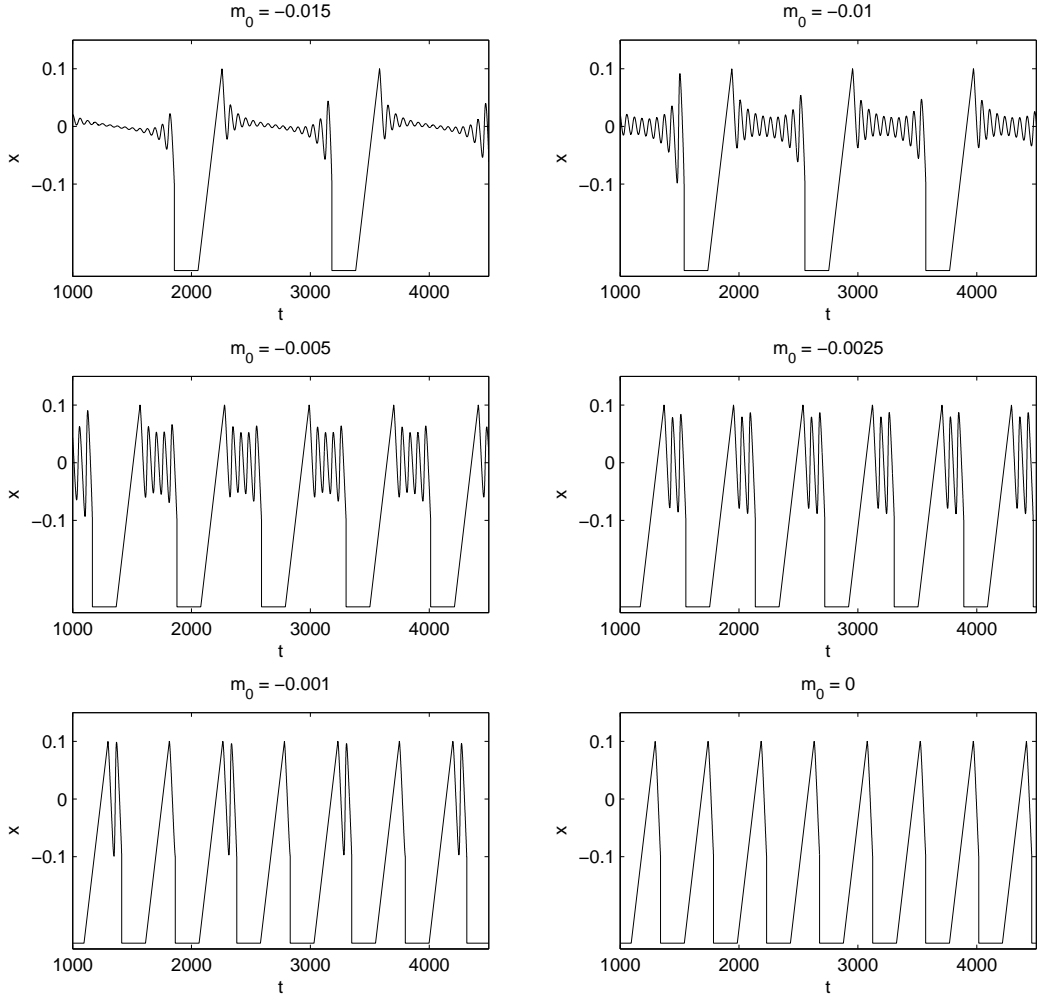}
	\caption{\label{fig:fig10}MMOs generated by the dynamical system with local dynamics \eqref{eq:BKW} (parameters $\epsilon=0.01$ and $\mu=0.006$) and global dynamics $m(z)=0.1z+m_0$ with $k_1=\sqrt\epsilon$. The transition from local to global dynamics has been applied when $x<-\sqrt\epsilon$. Trajectories have been started at $(x,y,z)=(\sqrt\epsilon,\epsilon,0.15)$.} 	
\end{figure}

\textit{Remark:} In the following, we are going to visualize the LAOs in a time series for the dynamical system defined by $m_{12}$ and $m_{21}$ by inserting a large amplitude oscillation at fixed amplitude whenever the map $m_{21}$ is applied.\\ 

For the folded node \eqref{eq:BKW} we fix the parameters $(\epsilon,\mu)=(0.01,0.006)$. Since $2k+1<\frac1\mu<2k+3$ for $k=82$ we know from folded node theory that there will be $k+2=84$ canards \cite{BronsKrupaWechselberger,Desrochesetal}. The theory also predicts that the maximum number of small oscillations is $k+1$. In Figure \ref{fig:fig10} we varied the global return mechanism to demonstrate that we can systematically reach sectors near a folded node with a sub-maximal number of oscillations and different MMO signatures. The global return map is chosen as the lowest order approximating linear map with $m(z)=0.1z+m_0$ and $A=0$, $b=0$. The parameter $m_0$ is viewed as the main bifurcation parameter and controls the entry of trajectories to the different folded node rotation sectors. We find the following MMO signatures:

\begin{center}
\begin{tabular}{|l|c|c|c|c|c|c|}  
\hline
$m_0=$ & $-0.015$ & $-0.01$ & $-0.005$ & $-0.0025$ & $-0.001$ & $0.0$ \\ 
$L^s=$ & $1^{14}$ & $1^9$   & $ 1^4$   & $1^2$     & $2^1$    & $1^0$ \\
\hline
\end{tabular}
\end{center}

Maximal MMO signatures can also be obtained but many of the small oscillations will be at an exponentially small scale due to the contraction towards the weak canard \cite{Desrochesetal}. Observe that we have efficiently de-coupled the local parameter dynamics from the global parameter dynamics.\\

\begin{figure}[htbp]
	\centering
		\includegraphics[width=1\textwidth]{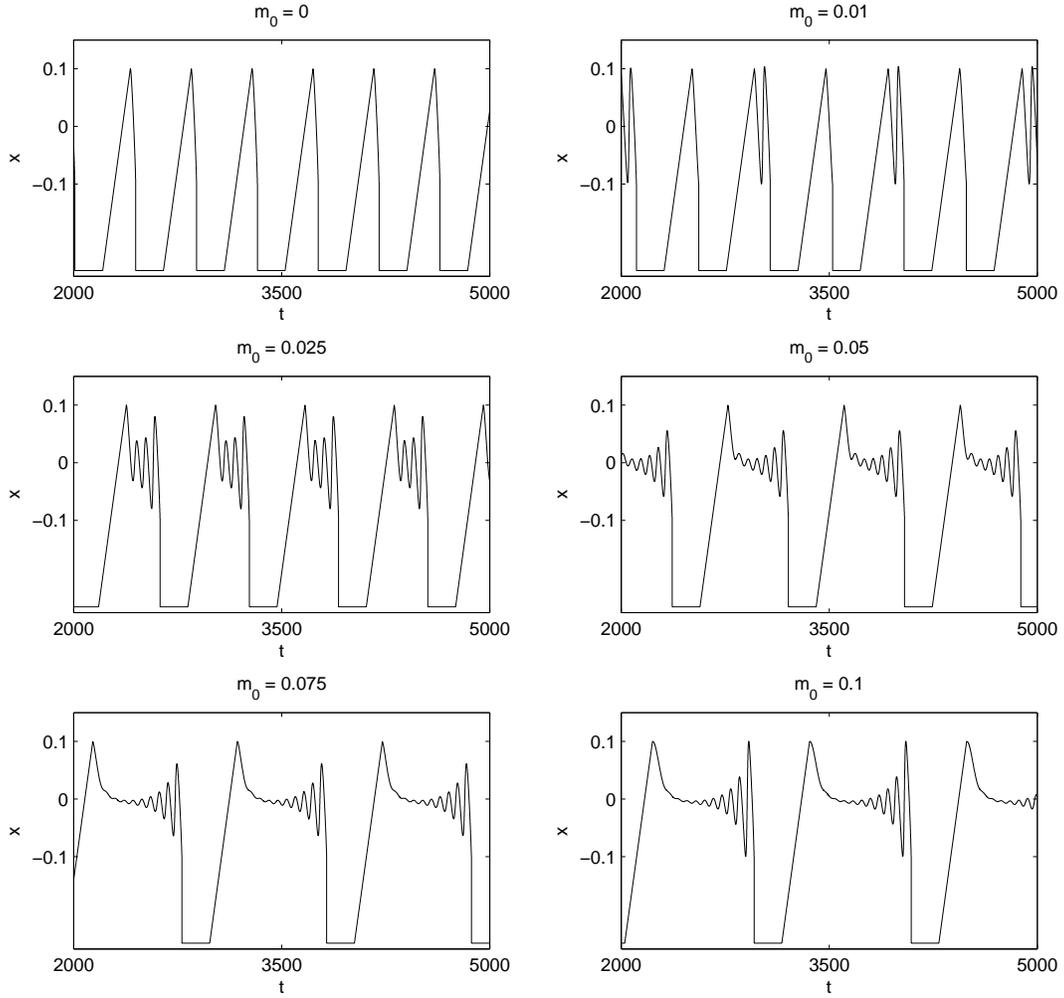}
	\caption{\label{fig:fig11}MMOs generated by the dynamical system with local dynamics \eqref{eq:Guck_SH} (parameters $(\epsilon,\nu,a,b,c)=(0.01,0.015,0.5, -1, 1)$) and global dynamics $m(z)=0.1z+m_0$ with $k_1=\sqrt\epsilon$. The transition from local to global dynamics has been applied when $x<-\sqrt\epsilon$. Trajectories have been started at $(x,y,z)=(\sqrt\epsilon,\epsilon,0.15)$.} 	
\end{figure}

The second simulation focuses on the singular Hopf normal form \eqref{eq:Guck_SH}. We fix the parameters $(\epsilon,\nu,a,b,c)=(0.01,0.01,0.5, -1, 1)$. For the global return map we again consider the singular limit of $m_{21}$ with parameter $m_2=0$, $m_1=0.1$ and primary bifurcation parameter $m_0$. It is easy to check that the equilibrium near the fold curve is located at $q=(x_{eq},y_{eq},z_{eq})\approx (-6.63729\times 10^{-3},4.40537\times 10^{-5},-6.63729\times 10^{-3})$. The equilibrium is a saddle-focus with one-dimensional stable and two-dimensional unstable manifold. We are in the regime where the SAOs are generated/amplified via the singular Hopf mechanism. Figure \ref{fig:fig11} shows the typical SAOs with increasing amplitude as they approach $W^u(q)$. We find the following MMO signatures:

\begin{center}
\begin{tabular}{|l|c|c|c|c|c|c|}  
\hline
$m_0=$ & $0.0$ & $0.01$ & $0.025$ & $0.05$ & $0.075$ & $0.1$ \\ 
$L^s=$ & $1^0$ & $2^1$   & $ 1^3$   & $1^6$     & $1^9$    & $1^9$ \\
\hline
\end{tabular}
\end{center}

It is important to note that the number of SAOs for the singular Hopf bifurcation is not only influenced by the folded node but also by the relative positions of the invariant manifolds of $q$ \cite{Desrochesetal}. In particular, Guckenheimer \cite{Guckenheimer7} points out that the one-dimensional stable manifold $W^s(q)$ seems to interact in an intricate way with MMO trajectories. Our decomposition approach is well-suited to investigate this dependency further once the local unfolding of the singular Hopf bifurcation is better understood \cite{GuckenheimerMeerkamp}. Also for the singular Hopf bifurcation we have been able to reproduce a typical MMO sequence without varying the local parameters. Extensive additional numerical simulation showed that it is difficult to produce periodic sequences of MMOs of the forms
\be
\label{eq:comp_MMOs}
L^s, \quad \text{with $L\gg 1, s\gg1$}\qquad \text{and} \qquad L_1^{s_1}L^{s_2}_2\cdots 
\ee
by varying further parameters in the map $m(z)$. These simulations confirm parts of the incomplete theory for MMOs in three dimensions \cite{KrupaPopovicKopell,KrupaPopovicKopellRotstein} which predict the limited number of MMO patterns for three time scale systems. Therefore we conjecture that higher-dimensional return maps are more likely to account for more complicated MMOs of the form \eqref{eq:comp_MMOs}.

\begin{figure}[htbp]
	\centering
		\includegraphics[width=0.98\textwidth]{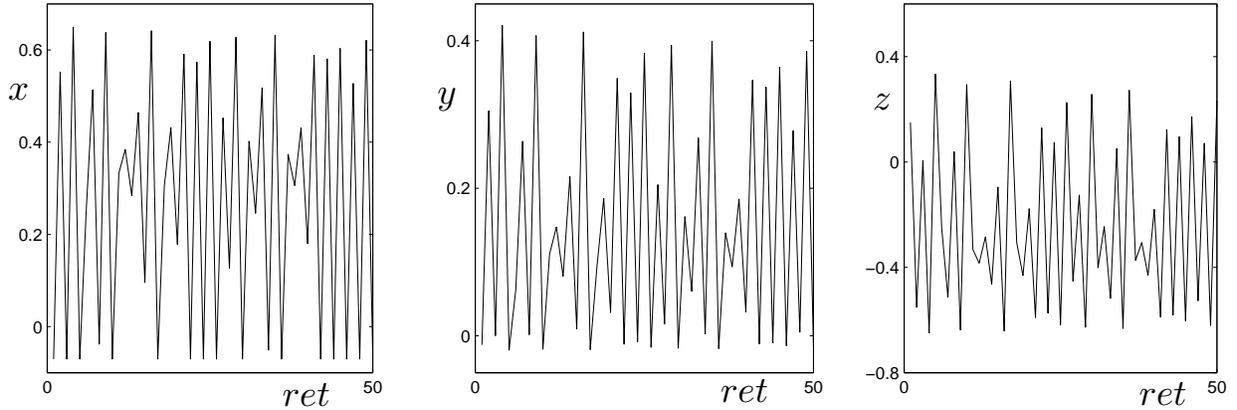}
	\caption{\label{fig:fig12}Coordinates $(x,y,z)$ before the global map is applied; the horizontal axis shows the return number i.e. $1^{st}$ application of the global map, $2^{nd}$ application, etc. The entire orbits are generated by the dynamical system with local dynamics \eqref{eq:Guck_SH} (parameters $(\epsilon,\nu,a,b,c)=(0.01,0.015,0.5, -1, 1)$) and global dynamics $m(z)=3z^2+0.2z-0.8m_0$ with $k_1=\sqrt\epsilon$. The transition from local to global dynamics has been applied when $x<-\sqrt\epsilon$. The times series of the returns shows typical chaotic non-periodic behavior.}	
\end{figure}

\begin{figure}[htbp]
	\centering
		\includegraphics[width=0.8\textwidth]{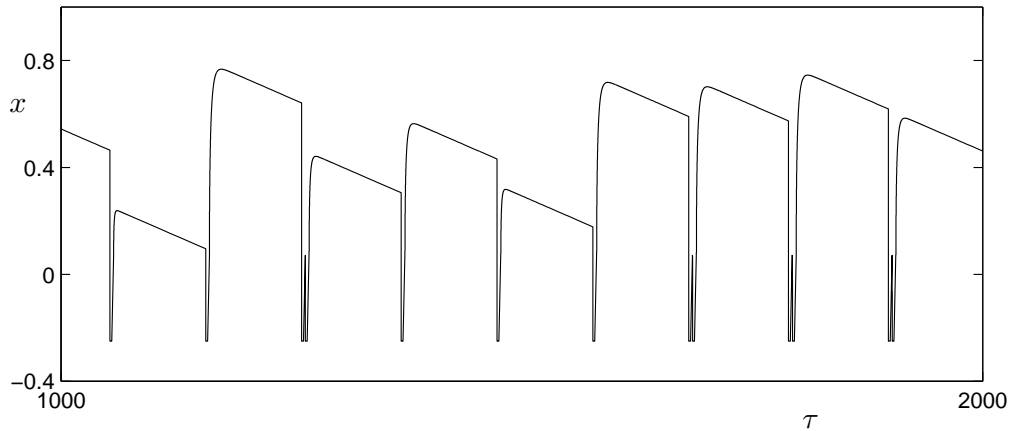}
	\caption{\label{fig:fig13} Subset of the time series in $(x,\tau)$ variables associated to the returns in Figure \ref{fig:fig12}. Irregular oscillations are observed with $1^0$ and $1^1$ components.}	
\end{figure}

Chaotic MMO signatures can be produced easily using a suitable quadratic map with $m_2\neq 0$. Figures \ref{fig:fig12} and \ref{fig:fig13} illustrate an orbit obtained from the dynamical system of the singular Hopf bifurcation with global returns generated by the map $m(z)=3z^2+0.2z-0.8$; the irregular behavior of the global returns in Figure \ref{fig:fig12} suggests that this orbit is chaotic. It is well-known that systems with two slow variables and one fast variable with S-shaped critical manifold can be chaotic \cite{MKKR,Haiduc1,GuckenheimerWechselbergerYoung,MedvedevYoo1}. Koper \cite{Koper} observed chaotic regions in parameter space in his original analysis of \eqref{eq:Koper}; transitions of MMOs to chaotic sequences can also be observed in many other models \cite{Desrochesetal}. As shown above, our model is also able to reproduces this aspect of typical MMO models. We conclude that our modeling approach reproduces the main dynamical features and decouples the global parameter dynamics from the local parameter dynamics.\\

\section{Brief Outlook}
\label{sec:outlook}

The strategy and methods we presented in this paper apply, in principle, to any system where the MMO mechanism can be decomposed into a local part that generates the SAOs and a global return map. For folded nodes and singular Hopf generated SAOs, the overall dimension can be arbitrary. Indeed, it has recently been shown by Wechselberger \cite{Wechselberger1} that the local theory in Appendices \ref{ssec:folded_nodes}-\ref{ssec:sing_Hopf} extends to systems with $m\geq 1$ fast and $n\geq 2$ slow variables. The technique is to use a center manifold reduction to get into the situation $(m,n)=(1,2)$. Then we can still compute singular return maps as we still have the three important one-dimensional curves that are analogous to $L^\mu$, $\gamma_s$ and $F^+$ in the Koper model. We have resolved the map for the Koper model in more detail using the drop curves $L^{a\pm}$. However, we could just compute $m_{a,+}\circ m_j$ and $m_{a,+}\circ m_f$ as single maps for another system or adapt the finer global decomposition to the fast-slow geometry of the problem. Moreover, it is very important to point out that a center manifold reduction has been already used in a four-dimensional system with MMOs generated by folded nodes and singular Hopf bifurcations \cite{Curtu}. Our methods apply verbatim to the resulting three-dimensional system obtained in \cite{Curtu}.

It is expected that the return maps for other systems can be more complicated. For example, just consider the situation for the Koper model but insert several non-trivial slow subsystem hyperbolic attractors on $C^{a-}$. Then the maps $m_j$ and $m_f$ may even have gaps since orbits can get trapped on persisting attractors on $C_\epsilon^{a-}$. Computing singular maps for several well-known MMO models \cite{Desrochesetal} and analyzing their structure is an interesting project but is beyond the scope of this paper.\\

\textbf{Acknowledgment:} I would like to thank two anonymous referees for valuable comments that helped to improve the focus and exposition of the paper.

\appendix

\section{Background Review}
\label{ssec:review}

\subsection{Fast-Slow Systems}
\label{ssec:fss}

We are only going to recall the basic definitions and results about fast-slow systems. There are several standard references that detail many parts of the theory \cite{Jones,KaperJonesPrimer,MisRoz,MKKR,Grasman,Desrochesetal,ArnoldEncy,KuehnBook}. A fast-slow system of ordinary differential equations (ODEs) is given by:
\be
\label{eq:basic1}
\begin{array}{rcrcl}
\epsilon \dot{x}&=&\epsilon\frac{dx}{d\tau}&=&f(x,y),\\
\dot{y}&=&\frac{dy}{d\tau}&=&g(x,y),\\
\end{array}
\ee 
where $x\in\R^m$ are fast variables, $y\in\R^n$ are slow variables and $0<\epsilon\ll1$ is a small parameter representing the ratio of time scales. Equation \eqref{eq:basic1} can be re-written by changing from the slow time scale $\tau$ to the fast time scale $t=\tau/\epsilon$
\be
\label{eq:basic2}
\begin{array}{lclcr}
x'&=&\frac{dx}{dt}=f(x,y),\\
y'&=&\frac{dy}{dt}=\epsilon g(x,y).\\
\end{array}
\ee
The singular limit $\epsilon\ra 0$ of \eqref{eq:basic2} yields the fast subsystem ODEs parametrized by the slow variables $y$. Setting $\epsilon\ra 0$ in \eqref{eq:basic1} gives a differential-algebraic equation (DAE), called the slow subsystem, on the critical manifold $C:=\{f(x,y)=0\}$. Concatenations of fast and slow subsystem trajectories are called candidates. 

A subset $S\subset C$ is called normally hyperbolic if the $m\times m$ total derivative matrix $(D_xf)(p)$ is hyperbolic. A normally hyperbolic subset $S$ is attracting if all eigenvalues of $(D_xf)(p)$ have negative real parts for $p\in S$; similarly $S$ is called repelling if all eigenvalues have positive real parts. On normally hyperbolic parts of $C$ the	implicit function theorem applies to $f(x,y)=0$ providing a map $h(y)=x$ so that $C$ can be expressed (locally) as a graph. Fenichel's Theorem \cite{Fenichel4,Jones,Tikhonov,WigginsIM} states that a compact normally hyperbolic submanifold $S=S_0$ of the critical manifold $C$ perturbs for $\epsilon>0$ sufficiently small, including stability and flow properties, to a slow manifold $S_\epsilon$.

A trajectory is called a maximal canard if it lies in the intersection of an attracting and a repelling slow manifold. Canards were first investigated by a group of French mathematicians \cite{BenoitCallotDienerDiener,Diener,Benoit1,Benoit4} using nonstandard analysis. Later also asymptotic \cite{Eckhaus,BaerErneux1,MKKR} and geometric \cite{DumortierRoussarie,KruSzm3,SzmolyanWechselberger1} methods have been developed to understand canard orbits.

\subsection{Folded Nodes}
\label{ssec:folded_nodes}

Normal hyperbolicity can fail in several ways. Here we briefly review the basic properties of two such situations \cite{Desrochesetal}. A non-degenerate fold point $p\in C$ is defined by requiring that $f(p) = 0$ and $(D_xf)(p)$ has rank $m - 1$ with left and right null vectors $w$ and $v$ so that $w \cdot [(D_{xx}f)(p)(v, v)] \neq 0$ and $w \cdot [(D_yf)(p)] \neq 0$. The set of fold points forms a manifold of codimension one in the $m$-dimensional critical manifold $C$. If $m=1$ and $n=2$ the fold points generically form a smooth curve that separates attracting and repelling sheets of the two-dimensional critical manifold $C$.

Two standard generating mechanisms for small oscillations near fold curves of the critical manifold will be considered in a normal form setup. Br{\o}ns, Krupa and Wechselberger \cite{SzmolyanWechselberger1,BronsKrupaWechselberger} consider a normal form 
\be
\label{eq:BKW}
\begin{array}{lcl}
\epsilon \dot{x}&=&y-x^2,\\
\dot{y}&=& -(\mu+1)x-z,\\
\dot{z}&=& \frac\mu2, \\
\end{array}
\ee 
where $x$ is the fast variable, $(y,z)$ are the slow variables and $\mu$ is a parameter. The critical manifold for \eqref{eq:BKW} is $C=\{y=x^2\}$ with a line of fold points $F=\{x=0,y=0\}$. $F$ decomposes the critical manifold $C=C^r\cup F\cup C^a$ where $C^r=C\cap \{x<0\}$ is repelling and $C^a=C\cap \{x>0\}$ is attracting. Differentiating $y=x^2$ implicitly with respect to $\tau$ gives $\dot{y}=2x\dot{x}$. Therefore the slow flow is
\be
\label{eq:BKW_sf}
\begin{array}{lcl}
\dot{x}&=& \frac{-(\mu+1)x-z}{2x},\\
\dot{z}&=& \frac\mu2 .\\
\end{array}
\ee   
Rescaling time by $\tau\mapsto 2x\tau$ reverses the direction of the flow on $C^r$ and yields the desingularized slow flow
\be
\label{eq:BKW_sf1}
\left(\begin{array}{c}
\dot{x}\\
\dot{z}\\
\end{array}\right)=\underbrace{\left(\begin{array}{cc}-(\mu+1) & -1 \\ \mu & 0\\\end{array}\right)}_{=:A_0} \left(\begin{array}{c}x \\ z\\ \end{array}\right)
\ee  
The desingularized slow flow has an equilibrium point at the origin $0=(0,0)\in F$ called a folded singularity. The eigenvalues $(\lambda_s,\lambda_w)=(-1,-\mu)$ of $A_0$ determine the type of the folded singularity. It is a folded saddle for $\mu<0$, a folded node for $\mu>0$ and a folded saddle-node of type II for $\mu=0$  \cite{SzmolyanWechselberger1,Desrochesetal}. We restrict to the folded node case and $\mu\in(0,1)$ here. Then $\lambda_s$ is associated to the strong eigendirection $\gamma_{s,0}$ and $\lambda_w$ is associated to the weak eigendirection $\gamma_{w,0}$. The extension of $\gamma_{s,0}$ ($\gamma_{w,0}$) under the slow flow is referred to as the strong (weak) singular canard. Trajectories in the funnel region bounded by $\gamma_{s,0}$ and $F$ can pass from $C^a$ to $C^r$; see also \cite{SzmolyanWechselberger1,BronsKrupaWechselberger}. 

The singular canards $\gamma_{0,s}$ and $\gamma_{0,w}$ perturb to maximal canards $\gamma_{\epsilon,s}$ and $\gamma_{\epsilon,w}$ that lie in the intersection of the two slow manifolds $C^a_\epsilon\cap C^r_\epsilon$ \cite{SzmolyanWechselberger1}. If $1/\mu\not\in \N$ then there are further maximal canards arising as intersections of $C^a_\epsilon\cap C^r_\epsilon$, called secondary canards \cite{Wechselberger2}. In particular, the attracting and repelling invariant manifolds twist around each other \cite{GuckenheimerHaiduc,Guckenheimer8}. The number of twists of a trajectory in the fold region can be predicted using its distance $\delta$ relative to the strong singular canard and by the value of $\mu$ \cite{BronsKrupaWechselberger}. We agree to the convention that $\delta>0$ indicates a trajectory entering the funnel region, $\delta=0$ describes the strong canard and for $\delta<0$ we are outside of the funnel. The twists can cause the SAOs of an MMO.

\subsection{Singular Hopf}
\label{ssec:sing_Hopf}

Note carefully that the normal form \eqref{eq:BKW} has no global equilibrium point for $\mu\in(0,1)$. However, in many applications a Hopf bifurcation occurs near the onset of MMOs \cite{Desrochesetal} which suggests to consider the possibility of a global equilibrium point passing through the folded node region. In particular, one has to add higher-order terms to the equation for $\dot{z}$ in \eqref{eq:BKW}. Augmenting these terms it is well-known that the global equilibrium can undergo a Hopf bifurcation at an $O(\epsilon)$-distance from the fold curve. This scenario is also been referred to as singular Hopf bifurcation \cite{Braaksma,Guckenheimer7} since the pair of complex conjugate eigenvalues involved in the Hopf bifurcation has a singular limit as $\epsilon\ra 0$ \cite{Braaksma}. Guckenheimer \cite{Guckenheimer7} derives the following normal form for a singular Hopf bifurcation
\be
\label{eq:Guck_SH}
\begin{array}{rcl}
\epsilon \dot{x}&=&y-x^2,\\
\dot{y}&=& z-x,\\
\dot{z}&=& -\nu-ax-by-cz, \\
\end{array}
\ee
where $x$ is a fast variable, $(y,z)$ are slow variables and $(\nu,a,b,c)$ are parameters. The key difference between \eqref{eq:BKW} and \eqref{eq:Guck_SH} is that we can find global equilibria $q=q(\nu,a,b,c)$ for \eqref{eq:Guck_SH}. They are determined by solving the equation
\be
\label{eq:ep_glob}
-\nu=(a+c)x+bx^2.
\ee  
If $\nu\approx 0$ then the equilibrium point is close to the folded singularity at the origin. The desingularized slow flow of \eqref{eq:Guck_SH} can be calculated similar to the folded node case. It can be shown \cite{Desrochesetal,KrupaWechselberger} that $q$ is only important for the local dynamics near $(0,0,0)$ if $\nu$ is smaller than $O(\epsilon^{1/2})$. The key difference between MMOs that pass near a global equilibrium is that the SAOs can also be influenced by the stable and unstable manifolds $W^s(q)$ and $W^u(q)$. Detailed visualizations of the situation can be found in \cite{Desrochesetal,DesrochesKrauskopfOsinga}. Results for the unfolding of \eqref{eq:Guck_SH} can be found in \cite{Guckenheimer7,GuckenheimerMeerkamp}. We are going to use the normal forms \eqref{eq:BKW} and \eqref{eq:Guck_SH} as ``black-box'' units for numerical simulation in Section \ref{sec:loc_glob}. 

\section{Error Analysis}
\label{sec:error}

\begin{figure}[htbp]
	\centering
		\includegraphics[width=0.8\textwidth]{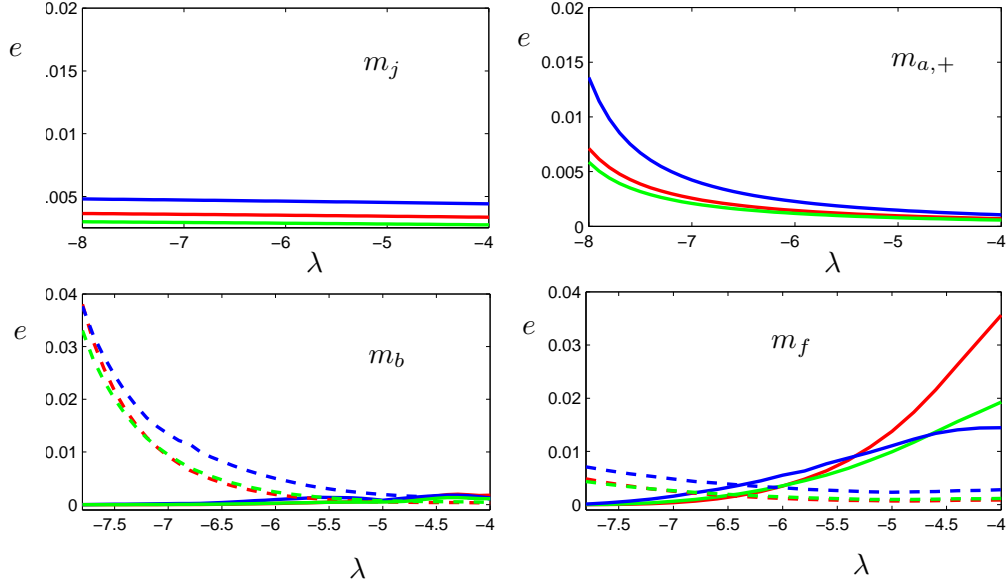}
	\caption{\label{fig:fig14}.Affine and quadratic fit approximation error $e$ of the different maps $m_{K}$ for $K\in\{j,(a,+),b,f\}$ given in \eqref{eq:error}. The horizontal axes are slices in $\lambda$-parameter space with $k=-10$. The error is measured in three different norms $L^1$ (red), $L^2$ (green) and $L^\I$ (blue). The dashed curves for $m_b$ and $m_f$ indicate the error for the lower branches and the solid curves the error for the upper branches of the maps. The domain for $m_{j}$ and $m_{a,+}$ has been chosen as $z\in(2\lambda-(4+k)-1,2\lambda-(4+k)+1)$ and the domain for $m_b$ and $m_f$ is the entire projection of the strong canard.}	
\end{figure}

We briefly analyze the error of our approximation for the maps $m_{K}$ for $K\in\{j,(a,+),b,f\}$ for $k=-10$ and $\lambda\in(\lambda_{FSN},\lambda_{nf})$. The numerical integration of trajectories was carried out with a standard stiff numerical integration method (\texttt{ode15s} in MatLab \cite{MatLab2010b_base}) with absolute error tolerance $10^{-8}$. The grid size $h$ for the domain of the maps $m_K$ was always chosen so that $h\leq 0.02$. The main question we have to address is whether at a given fixed set of parameters $(\lambda,k)$ there exist affine and quadratic approximations as postulated in Section \ref{sec:maps2}. Figure \ref{fig:fig14} shows the error of the fit to the postulated polynomial forms measured in three different norms 
\be
\label{eq:error}
\begin{array}{lcl}
e(L^1)&=&\int_{z_0}^{z_1} |m_K^{num}(z)-m_K^{fit}(z)|dz,\\
e(L^2)&=&\left(\int_{z_0}^{z_1} (m_K^{num}(z)-m_K^{fit}(z))^2dz\right)^{1/2},\\
e(L^\I)&=&\sup_{z\in[z_0,z_1]} |m_K^{num}(z)-m_K^{fit}(z)|,\\
\end{array} 
\ee   
where $m_K^{num}$ indicates the map obtained from numerical integration and $m^{fit}_K$ denotes the affine and quadratic fits. The integrals in \eqref{eq:error} have been evaluated from the discrete numerical integration data and the associated polynomials fits using a composite Simpson rule \cite{StoerBulirsch} which has error $\mathcal{O}(h^5)$ as $h\ra 0$. Figure \ref{fig:fig14} shows that the worst-case error for the proposed affine and quadratic maps due is at most on the order of $10^{-2}$ over the entire range of parameters; the numerical integration error $10^{-8}$ and the numerical quadrature error $h^5\leq (0.02)^5$ can be neglected here. Overall, the affine and quadratic approximations are certainly satisfactory to extract the basic MMO patterns.

A natural question is to ask what happens to the perturbations of $m_K$ when $\epsilon>0$. It is well-known from Fenichel theory that the error near normally hyperbolic segments of the critical manifold and in the fast subsystem is at most $\mathcal{O}(\epsilon)$ as $\epsilon\ra 0$. Near the fold points \cite{SzmolyanWechselberger} it has been proven that the error is at most $\mathcal{O}(\epsilon^{1/3})$ as $\epsilon\ra 0$. Therefore we have that $m_K(z)+\mathcal{O}(\epsilon^{1/3})$ represents a flow map for $0<\epsilon \ll1$. 

The numerical computations we present here can likely be made mathematically rigorous \cite{Haiduc1} using interval arithmetic and tools such as IntLab \cite{Rump}. The main reason for this conjecture is that rigorous numerical integration and quadrature are two standard situations in interval arithmetic \cite{Rump}. However, carrying out this rigorous proof is beyond the scope and goals of this paper.

\end{document}